\documentclass[11pt,reqno]{amsart}
\usepackage{amssymb,latexsym}
\usepackage{graphicx}

\begin{document}

\setcounter{page}{1}

\title[Fibonacci Sequence: Some Observations and A new Constant]{Convolutions induced Discrete Probability Distributions and a new Fibonacci Constant}

\author{Arulalan Rajan, Jamadagni, Vittal Rao}
\address{Centre for Electronics Design and Technology\\
                Indian Institute of Science\\
                Bangalore, India}             
\email{mrarul@cedt.iisc.ernet.in}


\author{Ashok Rao}
\address{Dept. of Electronics and Communication Engg\\
         Channabasaweswara Institute of Technology\\
         Gubbi, India}
\begin{abstract}
This paper proposes another constant that can be associated with Fibonacci sequence. In this work, we look at the probability distributions generated by the linear convolution of Fibonacci sequence with itself, and the linear convolution of symmetrized Fibonacci sequence with itself. We observe that for a distribution generated by the linear convolution of the standard Fibonacci sequence with itself, the variance converges to 8.4721359\ldots. Also, for a distribution generated by the linear convolution of symmetrized Fibonacci sequences, the variance converges in an average sense to 17.1942\ldots, which is approximately twice that we get with common Fibonacci sequence. 
\end{abstract}

\maketitle

\section{Introduction}
Fibonacci sequence is known to have some fascinating properties. The most classical of these is the convergence of the ratio $\frac{f[n]}{f[n-1]} $to the golden mean $\varphi$ \cite{kepler},\cite{knott}. Viswanath \cite{viswa}, proved that with certain randomness introduced in Fibonacci sequence,  the $n^{th}$ root of the absolute value of the $n^{th}$ term in the sequence converges,with a probability of 1 (i.e. with extremely rare exceptions, almost surely), to another constant 1.13198824 \ldots. In this paper, we propose yet another constant that can be associated with the Fibonacci sequence. This constant, corresponds to the limit variance of the distribution generated by the convolution of Fibonacci sequence with itself. The variance of such a distribution asymptotically approaches to 8.4721359 \ldots.\\
The paper is organized as follows:\\
Section \ref{const} of the paper discusses, in brief, some of the known constants associated with the Fibonacci sequence. In section \ref{fibdist} we propose the use of Fibonacci sequence convolution to generate discrete probability distributions. In section \ref{vari} we discuss the results followed by a few observations in section \ref{obs}, with conclusions in section \ref{conclude}.

\section{Some Constants Associated with Fibonacci Sequence}
\label{const}
It is well known that the classical Fibonacci sequence is generated by the recursion given in eq.\ref{fibonacci} below\cite{sloane}
\begin{equation}
\label{fibonacci}
f[n] = f[n-1]+ f[n-2] \hspace{0.3in} \textup{where} f[0] = 0; f[1] = 1.
\end{equation}
A closed form expression for generating Fibonacci sequence is given by the Binet's formula, as in Eq.\ref{binet}, below, that relates golden mean $\varphi$ to the sequence.
\begin{equation}
\label{binet}
f\left[n\right] = {{\varphi^n-(1-\varphi)^n} \over {\sqrt 5}}={{\varphi^n-(-1/\varphi)^{n}} \over {\sqrt 5}}
\end{equation}
where $\varphi = \frac{1+\sqrt5}{2}$. This is one of the most important and well known constants that has been associated with the Fibonacci sequence ever since Kepler showed that the golden ratio is the limit of the ratios of successive terms of the Fibonacci sequence.
A decade ago, Viswanath \cite{viswa} computed another constant that explained how fast the random Fibonacci sequences grow. In his work, he introduced randomness in the recurrence that generates Fibonacci sequence. Introducing randomness in Fibonacci sequence, we get a random Fibonacci sequence that is defined by the recurrence $f[n] = f[n-1]\pm f[n-2]$ with signs chosen as given below in Eq. \ref{rndfibo}, 
\begin{eqnarray}
\label{rndfibo}
 f[n]= \begin{cases} 
f[n-1]+f[n-2],  \mbox{with probability 0.5}; \\ 
f[n-1]-f[n-2],  \mbox{with probability 0.5}. 
\end{cases}
\end{eqnarray}
Viswanath showed that in the case of the above random Fibonacci sequence,
\begin{equation}
\sqrt[n]{|f[n]|} \to 1.13198824\dots \mbox{ with probability 1, as } n \to \infty.
\end{equation}
In this work, we propose yet another constant that is closely related to the Fibonacci sequence.
\section{Convolution of Fibonacci Sequence and Discrete Probability Distributions}
\label{fibdist}
\begin{figure}[h]
	\begin{center}
	\includegraphics[scale=0.5]{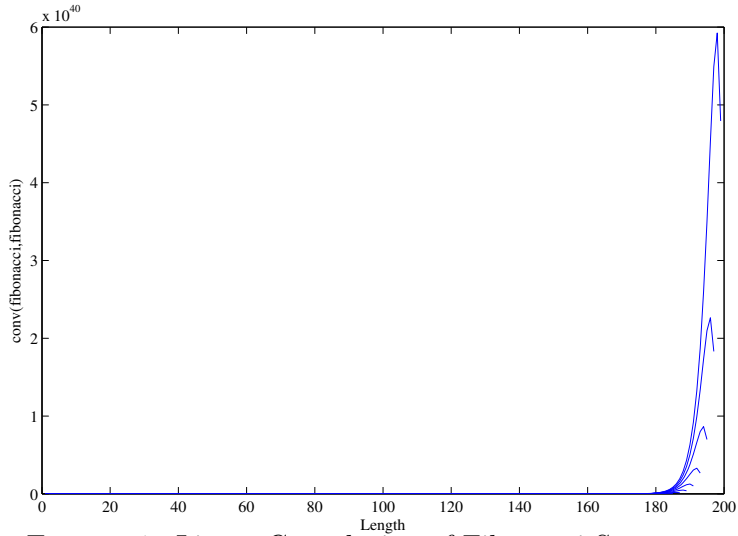}
	\end{center}
	\vspace{-0.28in}
	\caption{Linear Convolution of Fibonacci Sequences}
	\label{fiboconv}
\end{figure}

\begin{figure}[h]
	\begin{center}
	\includegraphics[scale=0.5]{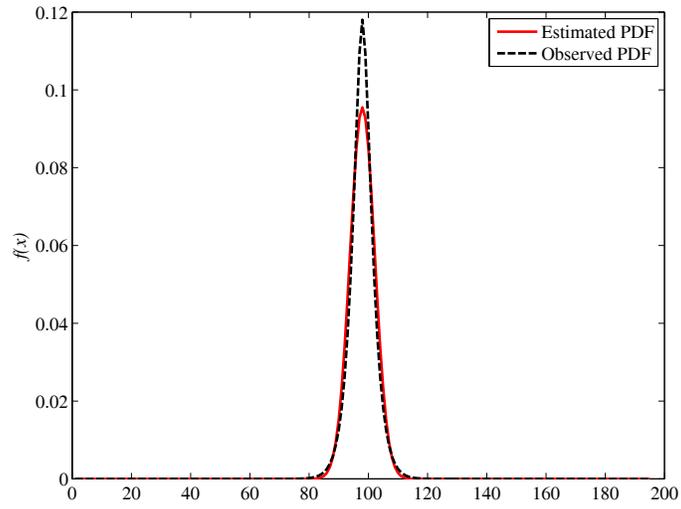}
	\end{center}
	\vspace{-0.26in}
	\caption{Fibonacci sequence (symmetry employed) convolution and narrow Gaussian distribution}
	\label{Fibocomp}
\end{figure}

 The motivation for this work comes from attempts to exploit the discrete nature of integer sequences for generating discrete probability distributions. In \cite{rajan}, the authors have looked at slow growing sequences and their convolutions, that approximates a Gaussian distribution with a mean squared error of about $10^{-8}$ or even less. In this section, we look at the discrete probability distributions generated by  a single convolution of Fibonacci sequences.\\
 Let $x_1[n]$ and $x_2[n]$ be two discrete sequences. Then the linear convolution of the two sequences is defined as, in Eq.(\ref{convd}),below. 
  \begin{equation}
  \label{convd}
  y[n] =  \sum_{k=-\infty}^{\infty} x_1[k]x_2[n-k] 
    \end {equation}    
 If $x_1[n]$ and $x_2[n]$ are two finite length sequences of length $L$ and $M$ respectively, then the length of $y[n]$ is $L+M-1$. If the two sequences to be convolved, have the same length $L$, then the length of $y[n]$ is $2L-1$. In this work, both $x_1[n]$ and $x_2[n]$ are finite length Fibonacci sequences of the same length, $L$.
  The convolution result is taken as the profile of the discrete probability distribution. The set of indices, ${n}$, namely ${1,2,3, \ldots, L+M-1}$, of $y[n] = x[n]*x[n]$, is considered the set of values that a discrete random variable $X$ can take. The probability of $X = n$ is defined by 
 \begin{equation}
       P(X = n) = \frac {y[n]}{\textup{sum}(y[n])},\hspace{0.1in}\textup{where,}\hspace{0.1in}\textup{sum}(y[n])=\sum_{n=1}^{L+M-1}y[n]
 \end{equation}
  Figure \ref{fiboconv} illustrates the convolution profiles obtained for different lengths of the classical Fibonacci sequence. Figure \ref{Fibocomp} gives the convolution of two Fibonacci sequences of the same length $L$, with symmetry employed at $L/2$. In this case, we generate Fibonacci sequence of length $L/2$ and then symmetrically extend it to length $L$. Figure \ref{Fibocomp} also compares the convolution result with an estimated normal probability density function with the same variance as the convolution result.
\section{Results}
\label{vari}
\begin{figure}[h]
	\begin{center}
	\includegraphics[scale=0.5]{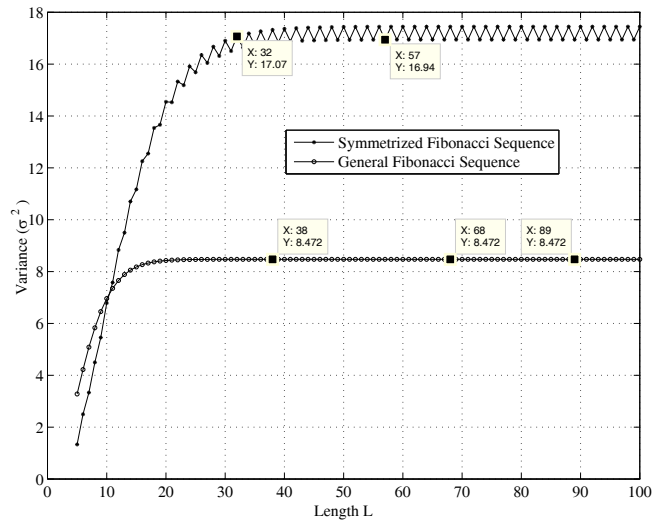}
	\end{center}
	\vspace{-0.2in}
	\caption{Variation of Variance with Length of Fibonacci Sequence}
	\label{FibVar}
\end{figure}

In the previous section, we looked at the linear convolution of Fibonacci sequence and the discrete distributions generated. We find that, increasing the length of Fibonacci Sequence and employing symmetry results in a profile that is similar to a very narrow Gaussian profile. For these distributions, we plot the variance and standard deviation as functions of the length of the sequence. From figure \ref{FibVar}, we observe that the variance of the distribution smoothly converge to a constant. This is due to the asymptotic and rapid growth of the Fibonacci sequence for large values of $n$. We also find from figure \ref{FibVar}, that the convergence, in case of Fibonacci sequence with symmetry employed, is only in the average sense.
\begin{figure}[h]
	\begin{center}
	\includegraphics[scale=0.5]{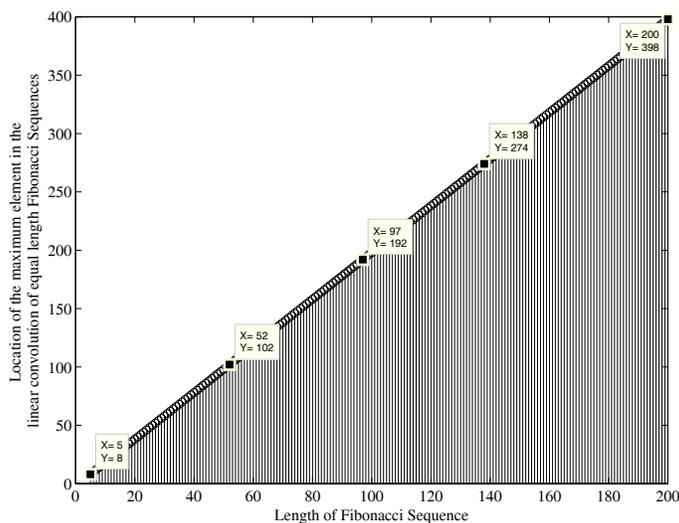}
	\end{center}
	\vspace{-0.2in}
	\caption{Plot of Index of Absolute Maximum for convolution of Fibonacci Sequences}
	\label{Fibmax}
\end{figure}
\section{Observations}
\label{obs}
We give below, some of the interesting properties, associated with Fibonacci sequence, that we have observed. Though, we have established a few of them mathematically, we are currently working on the analysis of the rest of them.
\begin{itemize}
\item Observation 1. With the standard Fibonacci sequence itself being considered as a distribution function, the variance converges to 4.23606797750108.

\item Observation 2. When we convolve (linear) standard Fibonacci sequence with itself and take the resulting sequence as a distribution function, the variance converges to 8.47213595500216.

\item Observation 3. The value mentioned in observation (2) is twice the value mentioned in observation (1).

\item Observation 4. Let $S_1$  be the standard Fibonacci Sequence. Let $S_2$ be the sequence $S_1$ in the reverse order. Let $S_3$ be the linear convolution of $S_1$ and $S_2$. The variance of $S_3$ converges to 8.47213595500216.

\item Observation 5. Let $S_1$  be the standard Fibonacci Sequence. Let $S_2$ be the sequence $S_1$ in the reverse order. Let $S_3$ be the linear convolution of $S_1$ and $S_2$. Let $S_4$ be a sequence obtained by linearly convolving $S_3$ with $S_1$ or $S_2$. The variance in either cases saturates to  12.7081989582623.  This value is 3 times the value listed in  observation (1).

\item Observation 6: The maximum value in the sequence resulting from the linear convolution of two standard Fibonacci sequences of length $L$ occurs at $2L-2$.

\item Observation 7:  Let $S_1$ be a standard Fibonacci, $L$ length sequence. Let $S_2$ be the sequence that is symmetrically extended version of $S_1$. $S_2$ has a length $2L$.  Convolving $S_2$ with itself yields another sequence $S_3$, of length $4L-1$. Considering $S_3$ as a distribution, the variance of $S_3$ converges in an average sense to 17.19423665579735. The swing is between 17.4442399455347 and 16.9442333660600.
\end{itemize}

From figure \ref{fiboconv}, we find that the sequence, resulting from the linear convolution of two increasing Fibonacci sequences of the same length, has its maximum at the extreme. From figure \ref{Fibmax}, we find that the maximum is located at $2L-2$. Mathematically, one can prove this by many ways. However, we follow the approach of proof by contradiction.
For this, it is enough to show that\\ 
     (I) $y[2L-3] < y[2L-2]$, \hspace{0.4in} (II) $y[2L-1] < y[2L-2]$ \\
where $y[2L-1]$ is the last element of the sequence $y[n]$ of length $2L-1$. This is because the linear convolution of two monotonically increasing functions will have only one absolute maximum. \\
\begin{proof}
First we prove $y[2L-1] < y[2L-2]$:
\begin{eqnarray}
\label{Y2L1}
y[2L-1]=&(f[L])^2,\\ 
\Rightarrow =&(f[L-1]+f[L-2])^2\\
\label{Y2L-1}
y[2L-1] =&(f[L-1])^2 + (f[L-2])^2 + 2f[L-1]f[L-2]
\end{eqnarray}
Now, we look at $y[2L-2]$.
\begin{eqnarray}
\label{Y2L2}
y[2L-2] &=& f[L-1]f[L] + f[L-1]f[L]\\
        &=& 2f[L-1]f[L]
\end{eqnarray}
We need to prove that Eq.(\ref{Y2L-1}) is less than Eq.(\ref{Y2L2})\newline
Let us assume that $y[2L-1] \geq y[2L-2]$ i.e.,
\begin{equation}
\label{inequal}
\footnotesize{(f[L-1])^2 + (f[L-2]))^2 + 2f[L-1]f[L-2] \geq 2f[L-1]f[L]}
\end{equation}
\begin{eqnarray}
\geq 2f[L-1]\bullet(f[L-1] + f[L-2]) \\
\geq 2f[L-1]^2 + 2f[L-1]f[L-2]
\end{eqnarray}
\begin{equation}
\label{in1}
\footnotesize{(f[L-1])^2 + (f[L-2])^2 + 2f[L-1]f[L-2]\geq 2(f[L-1])^2 + 2f[L-1]f[L-2]}
\end{equation}
From (\ref{in1}) we find that, 
\begin{equation}
\label{contra1}
   (f[L-2])^2 \geq (f[L-1])^2 
\end{equation}
    This implies that $f[L-2] > f[L-1]$ which is a contradiction arising due to our assumption that $y[2L-1] \geq y[2L-2]$. Therefore,
\begin{equation}
y[2L-1] < y[2L-2] \hspace{0.2in}\forall \hspace{0.1in} L>3
\end{equation}
We now look at the first inequality $y[2L-3] < y[2L-2]$.Let us assume that $y[2L-3] \geq y[2L-2]$ i.e.,
\begin{eqnarray}
y[2L-3] = 2f[L-2]f[L] + (f[L-1])^2\\
\label{inequal}
\Rightarrow 2f[L-2]f[L] + (f[L-1])^2  \geq 2f[L-1]f[L]
\end{eqnarray} 
Substituting for $f[L] = f[L-1] + f[L-2]$, we get,
\begin{equation}
\footnotesize{2(f[L-2])^2 + 2f[L-1]f[L-2] + (f[L-1])^2 \geq 2(f[L-1])^2 + 2f[L-1]f[L-2]}\\
\end{equation}
\begin{eqnarray}
\label{in2} 
2(f[L-2])^2 \geq (f[L-1])^2 \\ 
\label{contra2} 
\Rightarrow \sqrt2 \geq \frac {f[L-1]}{f[L-2]}
\end{eqnarray} 
This is a contradiction, as the ratio of $\frac {f[L-1]}{f[L-2]}$ approaches the golden ratio $\varphi$ which is greater than $\sqrt2$. Thus we find that on convolving Fibonacci sequence of length $L$, ($L \geq 3$), with itself, results in a sequence $y[n]$, that has its maximum at $n = 2L-2$.
\end{proof}
From the previous section \ref{vari} we note that the variance of a discrete distribution generated by the linear convolution of Fibonacci sequence with itself, saturates to a constant of value 8.4721359\ldots.
For a distribution generated by the linear convolution of symmetrized Fibonacci sequences, the variance saturates in an average sense to 17.1942\ldots, which is approximately two times 8.4721359\ldots. 
\section{Conclusion}
\label{conclude}
This work proposes another constant 8.4721359\ldots that can be associated with Fibonacci sequence. The constant corresponds to the value to which the variance, of the discrete distribution generated by linear convolution of Fibonacci sequences of the same length, converges to. Further, it is interesting to observe that on performing linear convolution of a Fibonacci sequence of length $L$ with itself, the maximum value always occurs at $2L-2$.

\medskip
\end{document}